\documentclass[12pt,fleqn]{amsart}

\usepackage{latexsym}
\usepackage{amssymb}
\usepackage{amsmath}
\usepackage{amsthm}
\usepackage{amsaddr}
\usepackage{amscd}
\usepackage{enumerate}
\usepackage{amsfonts}
 \usepackage{anysize}
\usepackage{color, soul}

\newtheorem{theorem}{Theorem}[section]
\newtheorem{lemma}[theorem]{Lemma}
\newtheorem{corollary}[theorem]{Corollary}
\newtheorem{proposition}[theorem]{Proposition}

\theoremstyle{definition}
\newtheorem{definition}[theorem]{Definition}
\newtheorem{remark}[theorem]{Remark}
\newtheorem{example}[theorem]{Example}

\numberwithin{equation}{section}
\newcommand{\N}{{\mathbb N}}

\newcommand{\n}{\noindent}
\newcommand{\vr}{\varepsilon}

\DeclareMathOperator{\cov}{Cov}

\marginsize{2cm}{2cm}{2cm}{2cm}

\title[Gromov-Hausdorff compactness theorem]{Gromov-Hausdorff
compactness theorem for non-Archimedean Fuzzy Metric 
Spaces}

 \author{{\sc Sergio Macario},\ \ \   {\sc Manuel Sanchis.}}
  
  \address{Institut Universitari de Matem\`{a}tiques i Aplicacions de
Castell\'{o}
 (IMAC), Universitat Jaume I, Campus del Riu Sec. s/n, 12071 Castell\'{o}
(Spain)}
 \email{macario@uji.es, sanchis@mat.uji.es}

 \date{}
  \thanks{Corresponding author: M. Sanchis, sanchis@uji.es
}
 \subjclass[2010]{}
 \keywords{Gromov´s compactness theorem, fuzzy metric space, Gromov-Hausdorff fuzzy metric}

\newcommand{\restringido}[1]{\,\vrule height5pt width.4pt depth10pt\,\lower
10pt\hbox{\scriptsize $#1$}}
  
\begin{document}
\maketitle

\begin{abstract}
The authors introduced in a previous paper the notion of fuzzy Gromov-Hausdorff distance between non-Archimedean compact fuzzy metric spaces, presenting a fuzzy version of the Gromov's completeness theorem. In this paper we present a fuzzy version of the Gromov's precompactness theorem which allows to deduce the classical theorem for compact metric spaces by means of the standard fuzzy metric associated to a metric space. This completes the  previous work.
\end{abstract}

\section{Introduction}

Gromov-Hausdorff convergence was introduced by Gromov in 
\cite{Gr} as a means of studying the convergence of metric spaces and it has turned out to be a powerful tool in areas like general topology,
geometry, functional analysis, biostatistical, etc... (\cite{an, boga, fu88, grpe, lee, memoli,pa, vik}). The interested reader might
consult the
survey \cite{pe} (see also \cite{BuBu}) on the Gromov-Hausdorff convergence of compact metric
spaces.

Since Gromov-Hausdorff metric has such a wide range of applications, it might be interesting to provide a Gromov-Hausdorff metric in the context of fuzzy metric spaces. The first step, which is to define a Hausdorff fuzzy metric, has been done by Rodr\'iguez-L\'opez and Romaguera in \cite{rorosa}. Thereafter, 
in \cite{MS:2015}, the authors introduced the so-called 
Gromov-Hausdorff fuzzy metric in the setting of compact non-Archimedean fuzzy metric spaces, which allowed them to show a fuzzy version of Gromov's completeness theorem.

In this paper we complete the previous work, giving a fuzzy version of the Gromov's precompactness theorem which allows to deduce the classical theorem for metric spaces by means of the standard fuzzy metric (defined below) associated to a metric space.

It is worth mentioning that one of the most useful tools in the theory of Gromov-Hausdorff convergence is the so-called Gromov's compactness theorem, an analog of the Arzelà-Ascoli theorem in the framework of uniform convergence. Gromov proved this result in his famous paper \cite{Gr} in 1981, and it has many applications.  For example, in sympletic topology, where the compactness of moduli spaces plays a key role \cite{Gr2}. Other interesting applications include the theory of pseudo-holomorphic curves \cite{Hu}, and Riemannian manifolds \cite{lee2}. For  other settings where this theorem is important, see, among others, \cite{ab, ac, O, wxy, Y}. Therefore, a fuzzy version of Gromov's compactness theorem might  open the way  for future applications in the fuzzy context. 
%
%

The paper is organized as follows. Sec\-tion~2 is devoted to present the
basic notions and facts about fuzzy metrics.
Section~3 sets out the main results about the non-Archimedean Gromov-Hausdorff fuzzy distance proved in our earlier paper \cite{MS:2015}. Section 4 presents the main results in this paper leading to  the  Gromov-Hausdorff's precompactness theorem and its relationship with the case of metric spaces. The conclusions 
are laid in the  last section.

\section{Preliminaries and Basic Facts}

We start with some fundamental concepts of fuzzy metric spaces used in this paper. A fuzzy set $M$ in a set $X$ is a map $M:X\rightarrow [0,1]$. To work with fuzzy sets, a $t$-norm is usually required. According to \cite{SchwSkla60}, a continuous $t$-norm is a binary operation $ \ast:\mathopen[0,1\mathclose] \times \mathopen[0,1\mathclose]\rightarrow \mathopen[0,1\mathclose]$ that meets the following criteria: (i) $\ast$ is associative and commutative, (ii) $\ast$ is continuous, (iii) $a\ast1=a$ for every $a\in \mathopen[0,1\mathclose]$, and (iv) $a\ast b\leq c\ast d$ whenever $a\leq c$ and $b\leq d$, with $a,b,c,d\in \mathopen[0,1\mathclose]$.

Given $\ast$ and $\diamond$, two continuous $t$-norms, the relationship $\ast\leq\diamond$ stands for $a\ast b\leq a\diamond b$, for all $a,b\in[0,1]$.  Examples of continuous $t$-norms are the minimum $a\wedge b=\min \{a,b\}$, 
the usual product $a\cdot b$ and the Lukasiewicz $t$-norm given by $a\ast_L b=\max\{a+b-1,0\}$. The relationship between them is $\ast_L\leq \cdot\leq \wedge$. As a matter of fact, $\ast\leq \wedge$, for each continuous
$t$-norm $\ast$.

\begin{definition}{\rm (\cite{KraMi75})}\label{def1} A \emph{fuzzy metric} (in
the sense of Kramosil and Michalek) on a set  $X$ is a pair $(M,\ast)$ such that
$M$ is a fuzzy set in $X\times X\times
 \mathopen[0,\infty\mathclose)$ and $\ast$ is a continuous  t-norm satisfying
for all $x,y,z\in X$ and $t,s>0$:

 \begin{enumerate}
  \item[(KM1)] $M(x,y,0)=0$;
 \item[(KM2)] $M(x,y,t)=1$ for all $t>0$ if and only if $x=y$;
 \item[(KM3)] $M(x,y,t)=M(y,x,t)$;
 \item[(KM4)] $M(x,y,t)\ast M(y,z,s)\leq M(x,z,t+s)$ and
 \item[(KM5)] $M(x,y,\cdot):\mathopen[0,\infty \mathclose)\rightarrow
 \mathopen[0,1\mathclose]$ is a left continuous function.
 \end{enumerate}
\end{definition}

 A triple $(X,M,\ast)$ such that $X$ is a set and $(M,\ast)$ is a
 fuzzy metric on $X$, will be called a a \emph{fuzzy metric space}. In $(X,M,\ast)$ we consider the topology $\tau_{M}$ where 
 a neighborhood basis for  $x\in X$  is formed by the sets:
 $$B_{M}(x,\varepsilon,t)=\{y\in X:M(x,y,t)>1-\varepsilon \}$$
 for all $\varepsilon \in \mathopen(0,1\mathclose)$ and $t>0$. These sets  are called open balls.

If the triangular inequality (KM4) of  Definition~\ref{def1} is replaced by 
\begin{enumerate}
\item[(NA) ] $M(x,z,\max\{t,s\})\geq M(x,y,t)\ast M(y,z,s)$
\end{enumerate} 
\noindent for all $x,y,z\in X$ and all $t,s>0$, then $(X,M,\ast)$ is called a
\emph{non-Archimedean fuzzy metric space} (see \cite{is}). It is routine 
 to show
that condition (NA) implies condition (KM4), that is, every non-Archimedean fuzzy
metric space is itself a fuzzy metric space. 
It is straightforward 
to verify that condition (NA) is equivalent to the two following
conditions: 
\begin{itemize}
\item[(NA1)] $M(x,z,t)\geq M(x,y,t)\ast M(y,z,t)$ and 
\item[(NA2)]
$M(x,y,\_)$ is nondecreasing for all $x,y\in X$.   
\end{itemize}

In the usual way, the following example shows how to relate \emph{a classic metric} with a fuzzy metric. 

 \begin{example}\label{cons1} {\rm(Compare \cite{GeVe94, GreRo00})}
 Let $(X,d)$ be a metric space. Define a fuzzy set
 $M_{d}$ in $X\times X\times \mathopen[0,\infty \mathclose)$ by
 $$
M_{d}(x,y,t)=\left\{
 \begin{array}{cl}
 \dfrac{t}{t+d(x,y)} & \text{ for all } x,y\in X \text{ and }
t>0\hbox{;}\medskip \\
 0 & \text{ for all } x,y\in X \text{ and } t=0\hbox{.}
 \end{array}
\right.
 $$

 It is apparent that   $(M_{d},\ast)$ is a fuzzy
metric on $X$ for all continuous $t$-norm $\ast$. It is called  the 
\emph{standard fuzzy metric} induced by $(X,d)$. Moreover, $(M\sb{d},\cdot)$ is non-Archimedean (for a more general situation, see 
\cite[Example~4.2.3]{alc}).   
On the contrary, $(M_d,\wedge)$ is a non-Archimedean fuzzy metric only when $d$ is an ultrametric (see \cite{gms}).

It can be easily shown that, being
 $\ast$ and $\diamond$ two continuous $t$-norms with $\diamond \leq \ast$, if $(X,M,\ast)$ is a non-Archimedean fuzzy metric space, then so is $(X,M,\diamond)$. Since the Lukasiewicz $t$-norm and the usual product $t$-norm satisfy $\ast_L\leq \cdot$, the above fact  implies that if $(X,d)$ is a metric space, then $(X,M_d,\ast_L)$ is also a non-Archimedean fuzzy metric space.
 \end{example}
 \begin{example}\label{ex:FMdiscrete}
 Let $(X,d)$ be a metric space. For every $\alpha\in[0,1)$ we can consider the fuzzy set 
 $$M_{\alpha 1}(x,y,0)=0\quad\text{and, for $t>0$,}\quad M_{\alpha 1}(x,y,t)=\begin{cases}
 \alpha &\text{if $t\leq d(x,y)$}\\
 1 & \text{if $t> d(x,y)$}\\
 \end{cases}$$
 Then $(X, M_{\alpha 1}, \ast)$ is a fuzzy metric space for every continuous $t$-norm $\ast$, although, in general, is not non-Archimedean since property (NA1) may fail. Nevertheless, when $X$ is a two-point space, $X=\{x_1,x_2\}$, then $(X,M_{\alpha 1},\ast)$ is a non-Archimedean fuzzy metric space.
 \end{example}

Another example of non-Archimedean fuzzy metric space is provided by the following notion.

\begin{definition}[\cite{grero04}]
A fuzzy metric $(M,\ast)$ on a set  $X$ is said to be stationary if $M$ does not depend on $t$, i.e. if for each $x, y \in  X$, the function $M_{(x,y)}(t) = M(x, y,t)$ is constant.
\end{definition}
It is apparent that any fuzzy metric space $(X,M,\ast)$ where $(M,\ast)$ is stationary is a non-Archimedean fuzzy metric space, since properties (KM4) and (NA) are then equivalent.

%
%
%




\section{The non-Archimedean Gromov-Hausdorff fuzzy distance}\label{sec3}

The non-Archimedean Gromov-Hausdorff fuzzy metric was introduced by the authors in \cite{MS:2015}, where they studied the convergence of Cauchy sequences of compact fuzzy metric spaces. Most of the concepts and results in this section can be found in that reference. We include it here for completeness. Next, we summarize the properties of such fuzzy metric. Our approach follows the pattern given in \cite{rorosa}. The first step is to define the notion of the Hausdorff fuzzy metric for compact subsets.

Let $A$ be a (nonempty) subset of a fuzzy
metric space $(X,M,\ast)$.  For $x\in X$ and $t>0$, let
$M(x,A,t)=\sup\{M(x,a,t)\, :\, a\in A\}$  (see
\cite[De\-fi\-ni\-tion~2.4]{vee}). 


\begin{definition}
Let $(X,M,\ast)$ be a fuzzy metric space. If $A$ and $B$ are (nonempty) compact subsets of 
$X$, then $$H\sb{M}(A,B,t)=\min \{\inf_{a\in A} M(a,B,t), \inf_{b\in B} M(A,b,t)\}$$ for all $t\geq 0$. 
\end{definition}
\begin{remark}\label{rem:Hm}
The definition above is a particular case of a more general definition for noncompact subsets (see \cite{rorosa}).  Notice that if  $A$ and $B$ are compact subsets of $X$ and  
$H_M(A,B,t)>1-\varepsilon$, the following conditions hold: 
\begin{enumerate}
    \item Given $a\in A$, there exists $b\in B$ with $M(a,b,t)>1-\varepsilon$. 
    \item Given $b\in B$, there exists $a\in A$ with $M(a,b,t)>1-\varepsilon$.
\end{enumerate}
Conversely, if, for some $t>0$ and some $0<\varepsilon<1$, conditions (1) and (2) are satisfied, then $H_M(A,B,t)\geq 1-\varepsilon$.
\end{remark}


 Given two
fuzzy metric spaces $(X, M_X, \ast)$ and $(Y, M_Y, \ast)$, 
a fuzzy metric $M$ on the disjoint union $X\sqcup Y$ 
is said to be \textit{admissible} if restricts to the given metric on $X$ and
$Y$,  respectively. In \cite{MS:2015}, the authors used this concept to reformulate the non-Archimedean  Gromov-Hausdorff fuzzy distance $(M_{GH},\ast)$ as follows.

\begin{definition}{\cite[Theorem 3.4]{MS:2015}}\label{thm:gh-admiss}
Given two non-Archimedean fuzzy metric spaces $(X, M_X, \ast)$ and
$(Y, M_Y, \ast)$, then
\[
M_{GH}(X,Y,t)=\sup \{H_M(X,Y,t) : M\in\mathcal{A}(X\sqcup Y)\}\quad (t\geq 0)
\]

\n where $\mathcal{A}(X\sqcup Y)$ stands for the set of all admissible non-Archimedean fuzzy metrics on $X\sqcup
Y$. 
\end{definition}
%
%
The non-Archimedean Gromov-Hausdorff
fuzzy metric 
satisfies conditions (i), (iii), (NA) and (v) in definition of a
non-Archimedean fuzzy metric. As pointed out in \cite{MS:2015},
condition (ii) may not be true.  However,
in the case of compact non-Archimedean fuzzy
metric spaces, condition (ii) is satisfied \emph{modulo isometry} (\cite[Theorem 3.11]{MS:2015}). 
Let $\mathcal{M}=\{(X\sb{j},M\sb{X\sb{j}},\ast) : j\in J\}$ be a set of compact
non-Archimedean fuzzy metric spaces, and let $\mathcal{I}(\mathcal{M})$ denote
the set of all isometry classes of elements of $\mathcal{M}$. Therefore,

\begin{theorem}{\cite[Theorem 3.14]{MS:2015}}
$(\mathcal{I}(\mathcal{M}),M_{GH},\ast)$ is a non-Archimedean fuzzy metric
space.  
\end{theorem}

\section{Gromov-Hausdorff fuzzy compactness theorem}

In \cite{MS:2015}, the authors proved that every Cauchy sequence of compact non-Archimedean fuzzy metric spaces is convergent in the non-Archimedean Gromov-Hausdorff fuzzy metric. Therefore, $(\mathcal{I}(\mathcal{M}),M_{GH},\ast)$ is a complete metric space.
Our goal in this paper is to complete that work by providing sufficient conditions to ensure that $(\mathcal{I}(\mathcal{M}),M_{GH},\ast)$ is a compact metric space. Being a complete space,  we only need  to study when it is precompact.

\begin{definition}
A sequence $\{(X\sb{n},M_{X_{n}},\ast)\}\sb{n\in\N}$ of non-Archimedean fuzzy
metric spaces converges to a non-Archimedean fuzzy metric space $(X,M_X,\ast)$
if $\lim\sb{n}M\sb{GH}(X_n,X,t)=1$ for all $t>0$.   
\end{definition}
\begin{definition}
A sequence $\{(X\sb{n},M_{X_{n}},\ast)\}\sb{n\in\N}$ of non-Archimedean fuzzy
metric spaces is said to be a Cauchy sequence if for every $t>0$ and $0<\varepsilon<1$ there exists $n_0\in\N$ such that $M\sb{GH}(X_n,X_m,t)>1-\varepsilon$ for all $n,m\geq n_0$.  
\end{definition}
\begin{remark}\label{rem:Cauchy}
It is worth remarking that a sequence $\{(X\sb{n},M_{X_{n}},\ast)\}\sb{n\in\N}$ of non-Archimedean fuzzy
metric spaces contains a Cauchy subsequence if, for every fixed $t>0$ and $0<\vr<1$, we can find a subsequence $\{(X\sb{n_k},M_{X_{n_k}},\ast)\}\sb{k\in\N}$ such that $M\sb{GH}(X_{n_j},X_{n_k},t)>1-\varepsilon$ for all $j,k\in\N$. Indeed, finding one of this for every $t=\frac{1}{n}$ and $\vr=\frac{1}{n}$, say $\{X_{n_k}^n\}_{k\in\N}$, provides the diagonal subsequence $\{X_{n_n}^n\}_{n\in\N}$  which satisfies that, for every  $t>0$ and $0<\vr<1$ and  $n,m\geq n_0$, where $\frac{1}{n_0}<\min\{\vr, t\}$,
$$M\sb{GH}(X_{n_n}^n,X_{m_m}^m,t)>M\sb{GH}\left(X_{n_n}^n,X_{m_m}^m,\frac{1}{n}\right)>1-\frac{1}{n}>1-\varepsilon.$$ 
\end{remark}
\begin{definition}
Let  $(X,M,\ast)$ be a fuzzy metric space. Given $0<\varepsilon<1$ and $t>0$,  we say that a set $A$ is a $(t,\varepsilon)$-net in $X$ if
for each $x\in X$ there exists $y\in A$ such that $M(x,y,t)>1-\varepsilon$.

\medskip

To simplify the notation, we will call it a $(t,\varepsilon\ast\varepsilon)$-net if $M(x,y,t)>(1-\varepsilon)\ast(1-\varepsilon)$ and so on.

The existence of a $(t,\varepsilon)$-net means that $X$ can be covered by balls  $B(y,\varepsilon,t)$ with centers $y$ in $A$.  
\end{definition}

\begin{definition} {\rm{\cite[Definition 1]{GreRo00}}}
A fuzzy  metric space $(X,M\sb{X},\ast)$ is called \emph{precompact} if for each $t>0$
and  
$0<\varepsilon<1$,  there is a finite subset $A$ of $X$, such that  $X =
\bigcup\sb{a\in A} B(a,\varepsilon,t)$. In this case, we say that $M_X$ is a precompact
fuzzy metric.  
This is equivalent to say that, for each  $t>0$
and  
$0<\varepsilon<1$,  there is a finite $(t,\varepsilon)$-net in $X$.
\end{definition}

Let $(X,M\sb{X},\ast)$ be a fuzzy metric space. Following \cite{MS:2015}, given $\delta\sb{1},
\delta\sb{2}$, with $0<\delta\sb{1}<1$ and $\delta_2>0$, we define  the ($\delta\sb{1}$,$\delta\sb{2}$)-\emph{cover number} of  $(X,M\sb{X},\ast)$ as 
\[
\cov (X,\delta\sb{1},\delta\sb{2})=\min \left \{ |C| :  X=\bigcup\sb{c\in
C}B(c,\delta\sb{1},\delta\sb{2}) \right\}.
\]
 
\n where $|C|$ stands for the cardinality of $C$. Notice that $\cov (X,\delta\sb{1},\delta\sb{2})$ is finite for all
$\delta\sb{1}, \delta\sb{2}$ if and only if $(X,M\sb{X},\ast)$  is precompact. 

Next result shows that, when $X$ and $Y$ are close enough,  we can find two finite nets in $X$ and $Y$ respectively, in such a way that the fuzzy distance between different elements of the nets in $X$ and $Y$ are close as well.

From now on, a fuzzy metric space is understood to be \emph{an Archimedean fuzzy metric space.}

\begin{proposition}\label{prop:necesaria}
Let $(X,M_X,\ast)$ and $(Y,M_Y,\ast)$ be two compact fuzzy metric spaces with $M_{GH}(X,Y,t)>1-\varepsilon$ for some $t>0$ and some $0<\varepsilon<1$. Then there exist two $(t,\varepsilon\ast\varepsilon\ast\varepsilon)$-net,   $\{x_i\}_{i=1}^{n}$  and $\{y_i\}_{i=1}^{n}$ in $X$ and $Y$, respectively, such that, for every $i,j=1,2\ldots, n$ we have that 
\begin{itemize}
\item[(a)] $M_X(x_i,x_j,t)\geq M_Y(y_i,y_j,t)\ast(1-\varepsilon)\ast(1-\varepsilon)$;
\item[(b)] $M_Y(y_i,y_j,t)\geq M_X(x_i,x_j,t)\ast(1-\varepsilon)\ast(1-\varepsilon)$.
\end{itemize}
\end{proposition}
\begin{proof}
Let  $t>0$ and $0<\varepsilon<1$ with $M_{GH}(X,Y,t)>1-\varepsilon$. By compacity of $X$ we can find a 
$(t,\varepsilon)$-net  in $X$, say $\{x_i\}_{i=1}^n$. Since $M_{GH}(X,Y,t)>1-\varepsilon$, we can fix an admissible fuzzy metric $M$ on $X\sqcup Y$ with $H_M(X,Y,t)>1-\vr$. 
Then, we can choose, for each $x_i$, an element $y_i \in Y$ with $M(x_i,y_i,t)>1-\varepsilon$ (see Remark~\ref{rem:Hm}).

\textsc{Claim:} $\{y_i\}_{i=1}^n$ is a $(t,\varepsilon\ast\varepsilon\ast\varepsilon)$-net in $Y$. To see this, notice that if $y\in Y$, then there exists an $x\in X$ with $M(x,y,t)>1-\varepsilon$. If we choose $x_i$ such that 
$M_X(x,x_i,t)>1-\varepsilon$, then  
\[
M_Y(y,y_i,t)= M(y,y_i,t)\geq M(x,y,t)\ast M(x,x_i,t)\ast M(x_i,y_i,t)>(1-\varepsilon)\ast(1-\varepsilon)\ast(1-\varepsilon).
\]

This prove the claim. Now (a) follows from the fact that, for every $i,j=1,2,\ldots, n$, we have 
\begin{align*}
M_X(x_i,x_j,t)&=M(x_i,x_j,t)\geq M(x_i,y_i,t)\ast M(y_i,y_j,t)\ast M(x_j,y_j,t)\\[2ex]
&>M_Y(y_i,y_j,t)\ast(1-\varepsilon)\ast(1-\varepsilon). 
\end{align*}
An argument similar to the previous one shows (b). 
\end{proof}

\begin{remark}
Notice that, being  $\{x_i\}_{i=1}^n$ a $(t,\varepsilon)$-net in $X$, it is a  $(t,\varepsilon\ast\varepsilon\ast\varepsilon)$-net as well.  
\end{remark}

\begin{definition}
Let $(X,M_X,\ast)$  be a fuzzy metric space. For each $t>0$, the $t$-diameter of $X$ is defined as
$${\rm diam}_t(X)=\inf\{ M_X(x,y,t)\, : \, x,y\in X\}.$$  
\end{definition}

If $(X,d)$ is a compact metric space and consider the fuzzy metric space $(X,M_d,\cdot$) with $M_d$ the standard metric 
associated to $(X,d)$, an easy check shows that, for all $t>0$, 
\[
{\rm diam}_t(X)=\dfrac{t}{t+{\rm diam} (X)}
\]

\noindent where ${\rm diam}(X)$ stands for the diameter of the metric space $(X,d)$. However, it is not difficult to find examples of compact fuzzy metric spaces $(X,M,\ast)$ such that ${\rm diam}_t(X)=0$ for infinite many values of $t$ (even for all $t>0$).

Assume that $(M,\ast)$ is a stationary fuzzy metric on $X$. Then the t-diameter does not depend on $t>0$, so it will be denoted by 
$${\rm diam}(X)=\inf\{ M_X(x,y)\, : \, x,y\in X\}.$$

\begin{lemma}\label{lem:Cs}
Let $(X,M_X,\ast)$ and $(Y,M_Y,\ast)$ be two  fuzzy metric spaces and assume that 
there exists a nondecreasing left-continuous function $C:[0,+\infty)\rightarrow [0,1]$, with 
$$C(s)\leq \min\{{\rm diam}_s(X),\, {\rm diam}_s(Y)\}\quad\text{ for all $s>0$.}$$
Then the fuzzy set $M$ on the disjoint union $X\sqcup Y$, defined as
\begin{enumerate}
\item $M(x,y,0)=0$, for all $x,y\in X\sqcup Y$;
\item $M(x,x',s)=M_X(x,x',s)$, for all $x, x'\in X$;
\item $M(y,y',s)=M_Y(y,y',s)$, for all $y, y'\in Y$;
\item $M(x,y,s)=M(y,x,s)=C(s)$, for all $x\in X$ and $y\in Y$,
\end{enumerate}
 is an admissible non-Archimedean fuzzy metric on $X\sqcup Y$.
\end{lemma}
\begin{proof}
It is apparent that properties (KM1),(KM2), (KM3), (KM5) and (NA2) are satisfied. Thus, we only need to show property (NA1). To see this, if $x,x'\in X$, $y,y'\in Y$ and $s>0$,  then  conditions (2) and (3) in the definition above imply that  we only need to consider the following cases:  

\[
\begin{aligned}
(a)\,\, M(x,x',s)&\geq C(s)\geq C(s)\ast C(s)=M(x,y,s)\ast M(x',y,s). \\[1ex]
(b) \,\, M(y,y',s)&\geq C(s)\geq C(s)\ast C(s)=M(x,y,s)\ast M(x,y's). \\[1ex]
(c) \,\, M(x,y,s)&=C(s)\geq M(x,x',s)\ast C(s)=M(x,x',s)\ast M(x',y,s). \\[1ex]
(d) \,\, M(x,y,s)&=C(s)\geq C(s)\ast M(y,y',s)=M(x,y',s)\ast M(y',y,s).\\[1ex]
\end{aligned}
\]
\end{proof}
\begin{remark}
It is worth noting that $C(s)=0$, for all $s>0$, satisfies  Lemma~\ref{lem:Cs}. In this case, the fuzzy metric $M$ defined in the previous lemma is always an admissible non-Archimedean fuzzy metric on $X\sqcup Y$.    
\end{remark}

\begin{example}
Let $(X,d_X)$,  $(Y,d_Y)$ be two compact metric spaces and consider the fuzzy metric spaces $(X,M_{d_X},\cdot)$ and 
$(Y,M_{d_Y},\cdot)$. If $K$ is a real number with $K\geq \max\{ {\rm diam(X)},\, {\rm diam}(Y)\}$, then the function defined as $C(s)=\frac{s}{s+K}$ for all $s>0$
satisfies the requested condition in the previous lemma. 
\end{example}

\begin{lemma}\label{lem:cond01}
Let $(X,M_X,\ast)$ and $(Y,M_Y,\ast)$ be two compact fuzzy metric spaces satisfying, for some $x,x'\in X$, $y,y'\in Y$, $t>0$ and $0<\varepsilon<1$, 
\begin{itemize}
\item[(a)] $M_X(x,x',t)\geq M_Y(y,y',t)\ast(1-\varepsilon)$,
\item[(b)] $M_Y(y,y',t)\geq M_X(x,x',t)\ast(1-\varepsilon)$.
\end{itemize}
Then, there exists $\delta>0$ such that, for every $t-\delta\leq s\leq t$, the properties (a) and (b) hold.
\end{lemma}
\begin{proof}
Take $0<\eta<\min\{M_X(x,x',t),M_Y(y,y',t)\}$ and, by the property (KM5) of the fuzzy metrics $M_X$ and $M_Y$, find $\delta>0$ such that, for all $t-\delta\leq s\leq t$, we have that
\begin{align*}
M_X(x,x',t)-M_X(x,x',s)&<\eta \quad\text{and}\quad M_Y(y,y',t)-M_X(y,y',s)<\eta.\\
\end{align*}
Then, for such $s$, we obtain that 
\begin{itemize}
\item[(a)] $M_X(x,x',s)\geq M_X(x,x',t)-\eta> M_Y(y,y',t)\ast(1-\varepsilon)-\eta\geq M_Y(y,y',s)\ast(1-\varepsilon)-\eta$,
\item[(b)] $M_Y(y,y',s)\geq M_Y(y,y',t)-\eta>M_X(x,x',t)\ast(1-\varepsilon)-\eta\geq M_X(x,x',s)\ast(1-\varepsilon)-\eta$,
\end{itemize}
for every such $\eta>0$. Thus, we have
\begin{itemize}
\item[(a)] $M_X(x,x',s)\geq M_Y(y,y',s)\ast(1-\varepsilon)$,
\item[(b)] $M_Y(y,y',s)\geq M_X(x,x',s)\ast(1-\varepsilon)$.
\end{itemize}
\end{proof}

Proposition~\ref{prop:necesaria} provides a necessary condition to get $M_{GH}(X,Y,t)>1-\varepsilon$. Next result shows that it is part of a sufficient condition as well.

\begin{proposition}\label{admisible}
Let $(X,M_X,\ast)$ and $(Y,M_Y,\ast)$ be two compact fuzzy metric spaces satisfying the following properties:
\begin{enumerate}
\item There exists a nondecreasing left-continuous function $C:[0,+\infty)\rightarrow [0,1]$, with $$C(s)\leq \min\{diam_s(X),\ diam_s(Y)\}\quad\text{ for all $s>0$.} $$
\item Suppose there exist $t>0$ and $0<\varepsilon<1$ such that we can find $n\in\mathbb{N}$ and  
two $(t,\varepsilon)$-nets,  $\{x_i\}_{i=1}^{n}$, $\{y_i\}_{i=1}^{n}$ in $X$ and $Y$ respectively, satisfying,  for all $s\geq t$ and for all $i,j=1,\ldots, n$,
\begin{itemize}
\item[(a)] $M_X(x_i,x_j,s)\geq  M_Y(y_i,y_j,s)\ast(1-\varepsilon)$,
\item[(b)] $M_Y(y_i,y_j,s)\geq M_X(x_i,x_j,s)\ast(1-\varepsilon)$.
\end{itemize}
\end{enumerate}
Then, there is an admissible non-Archimedean fuzzy metric $M$ on $X\sqcup Y$ with
$H_{M}(X,Y,t)>(1-\varepsilon)\ast(1-\varepsilon)$. Therefore, $M_{GH}(X,Y,t)>(1-\varepsilon)\ast(1-\varepsilon)$.
\end{proposition}
\begin{proof}
Take  $\{x_i\}_{i=1}^{n}$  and $\{y_i\}_{i=1}^{n}$  two $(t,\varepsilon)$-nets given by condition $(2)$.
Applying Lemma~\ref{lem:cond01}, we find $\delta>0$ such that, for $t-\delta\leq s\leq t$ and for all $i,j$, we also have:
\begin{itemize}
\item[(a)] $M_X(x_i,x_j,s)\geq M_Y(y_i,y_j,s)\ast(1-\varepsilon)$.
\item[(b)] $M_Y(y_i,y_j,s)\geq M_X(x_i,x_j,s)\ast(1-\varepsilon)$.
\end{itemize}

Let $M_{\delta}:X\times X\times [0,+\infty)\rightarrow [0,1]$ be defined as
\begin{enumerate}
\item $M_{\delta}(x,y,0)=M_{\delta}(y,x,0)=0$, for all $x\in X$ and $y\in Y$.
\item $M_{\delta}(x,x',s)=M_X(x,x',s)$, for all $s\geq0$ and $x,x'\in X$.
\item $M_{\delta}(y,y',s)=M_Y(y,y',s)$, for all $s\geq 0$ and $y,y'\in Y$.
\item For all $x\in X$ and $y\in Y$ and for all $s>0$,
\[
M_{\delta}(x,y,s)=M_{\delta}(y,x,s)=
\begin{cases}
C(s)\ast C(s)\ast (1-\varepsilon), &\text{if $0<s\leq t-\delta$},\\[2ex]
\text{\small $\displaystyle\max_{1\leq j\leq n}\Big\{ M_X(x,x_j,s)\ast M_Y(y,y_j,s)\Big\}\ast(1-\varepsilon)$},
&\text{if $s>t-\delta$}.\\[2ex]
\end{cases}
\]
\end{enumerate}
 In  order to prove that $M_{\delta}$ is a non-Archimedean fuzzy metric, we need  to check  properties (KM1), (KM2), (KM3), (KM5) and (NA) in  Definition~\ref{def1}. Properties (KM1), (KM2), (KM3)  are trivially satisfied. We show condition (NA) by verification of properties (NA1) and (NA2).
 
 It is apparent that property (NA2) is satisfied, since, for all $x\in X$,  $y\in Y$  and $j=1,\ldots,n$, $M_X(x,x_j,\_)$ and $M_Y(y,y_j,\_)$ are nondecreasing functions.
 
 For property (NA1) we need to consider two possibilities 
 depending on where the value $s>0$ is placed. Firstly, fix $s>t-\delta$ and consider 
 the following cases: 
 
 \textsc{Case 1}: $M_{\delta}(x,x',s)\geq M_{\delta}(x,y,s)\ast M_{\delta}(x',y,s)$, $x,x'\in X$, $y\in Y$. 
 For, choose $j$ and $k$ such that
\[
\begin{aligned}
 M_{\delta}(x,y,s) &= M_X(x,x_j,s)\ast M_Y(y,y_j,s)\ast(1-\varepsilon)  \text{ and }\\
 M_{\delta}(x',y,s) &=M_X(x',x_k,s)\ast M_Y(y,y_k,s)\ast(1-\varepsilon) 
 \end{aligned}
 \]
Then,
\[
\begin{split}
M(x,x',s)&=M_X(x,x',s)
\geq M_X(x,x_j,s)\ast M_X(x_j,x_k,s)\ast M_X(x_k,x',s)\\[2ex]
&\geq M_X(x,x_j,s)\ast M_Y(y_j,y_k,s)\ast (1-\varepsilon)\ast M_X(x_k,x',s)\\[2ex]
&\geq M_X(x,x_j,s)\ast M_Y(y,y_j,s)\ast M_Y(y,y_k,s)\ast(1-\varepsilon)\ast M_X(x_k,x',s)\\[2ex]
&=M_{\delta}(x,y,s)\ast M_Y(y,y_k,s)\ast M_X(x_k,x',s)\\[2ex]
&\geq  M_{\delta}(x,y,s)\ast  M_Y(y,y_k,s)\ast M_X(x_k,x',s)\ast (1-\varepsilon)\\[2ex]
&=M(x,y,s)\ast M(x',y,s)
\end{split}
\]

\textsc{Case 2}: $M_{\delta}(y,y',s)\geq M_{\delta}(x,y,s)\ast M_{\delta}(x,y',s)$, $x\in X$, $y,y'\in Y$. 
Analogously to the previous one but using property (b) instead of (a).

\textsc{Case 3}: $M_{\delta}(x,y,s)\geq M_{\delta}(x,x',s)\ast M_{\delta}(x',y,s)$, $x,x'\in X$, $y\in Y$. 

For, choose $1\leq k\leq n$ such that
$$M_{\delta}(x,y,s)=M_X(x,x_k,s)\ast M_Y(y,y_k,s)\ast(1-\varepsilon)$$
Then,
\[
\begin{split}
M_{\delta}(x,y,s)&\geq M_X(x,x_k,s)\ast M_Y(y,y_k,s)\ast(1-\varepsilon)\\[2ex]
&\geq  M_X(x,x',s)\ast M_X(x',x_k,s)\ast M_Y(y,y_k,s)\ast(1-\varepsilon)\\[2ex]
&=M_{\delta}(x,x',s)\ast M_{\delta}(x',y,s)
\end{split}
\]

\textsc{Case 4}: $M_{\delta}(x,y,s)\geq M_{\delta}(x,y',s)\ast M_{\delta}(y',y,s)$, $x\in X$, $y,y'\in Y$. Similarly to the previous case.

Notice that conditions (2) and (3) in the definition  of $M_{\delta}$ gives 
 the property (NA1) for the other possibilities. Therefore, property (NA1) is true for $s>t-\delta$. 
 
Otherwise, if $0<s\leq t-\delta$, the property can be proved as in Lemma~\ref{lem:Cs}, since $C(s)\ast C(s)\ast (1-\varepsilon)\leq C(s)\leq \min\{diam_s(X),\ diam_s(Y)\}$.


\medskip
 Condition (KM5) follows from the fact that $C(s)$ is a left-continuous and \cite[Proposition~2.6]{MS:2015}.

\medskip
Finally, we have to show that $H_M(X,Y,t)> (1-\varepsilon)\ast(1-\varepsilon)$. To see this, 
take any $x\in X$ and find $x_j$ in the $(t,\varepsilon)$-net in $X$ such that $M(x,x_j,t)=M_X(x,x_j,t)>1-\varepsilon$.
Then, we get 
\[
M(x,y_j,t)\geq M_X(x,x_j,t)\ast M_Y(y_j,y_j,t)\ast (1-\varepsilon) > (1-\varepsilon)\ast(1-\varepsilon).
\]
Similarly, given any $y\in Y$ we can find $y_j$ with $M(y,y_j,t)>(1-\varepsilon)$; so 
\[
M(y,x_j,t)\geq  M_X(x_j,x_j,t)\ast M_Y(y,y_j,t)\ast (1-\varepsilon) > (1-\varepsilon)\ast(1-\varepsilon).
\]
Then, $M_{GH}(X,Y,t)\geq H_M(X,Y,t)> (1-\varepsilon)\ast(1-\varepsilon)$ and we are done.
\end{proof}

Therefore, in order to show that two compact spaces are close in the Gromov-Hausdorff fuzzy metric, we need to find $(t,\varepsilon)$-nets satisfying condition (2) of the previous proposition. 

\medskip
We will work with continuous $t$-norms $\ast$ satisfying the following property: \begin{equation}\label{eq:distrib}
a-a\ast b\geq a\ast(1-b),\quad \text{for every $a,b\in[0,1]$.}
\tag{TN1}
\end{equation}
The product and the  Lukasiewicz $t$-norms are examples of such types of continuous $t$-norms. On the contrary, the minimum $t$-norm does not satisfy property \eqref{eq:distrib}. But it is worth mentioning that fuzzy metric spaces with the minimum $t$-norm $\wedge$ can fail to be non-Archimedean.

\begin{lemma}\label{lem:ab}
Let $(X,M_X,\ast)$ be a fuzzy metric space with a $t$-norm $\ast$ satisfying property \eqref{eq:distrib}. Let $a$ and $b$ be two numbers in $(0,1)$ and assume that $0<K<\min\{a,b\}<1$. Let  $0<\varepsilon<1$.  If $|a-b|<K\ast \varepsilon$, then
$a\geq b\ast(1-\varepsilon)$ and $b\geq a\ast(1-\varepsilon)$.
\end{lemma}
\begin{proof}
When $a=b$ the property is trivially true. We can assume, without loss of generality, that $a>b$. Then
$$a> b\geq b\ast(1-\varepsilon)$$
Besides, $|a-b|=a-b<K\ast\varepsilon$ and then
$$b>a-K\ast\varepsilon \geq a-a\ast\varepsilon\geq a\ast(1-\varepsilon)$$
\end{proof}
It is worth to remark that if $a\neq b$, then we get the strict inequalities
$$a>b\ast(1-\varepsilon)\text{ and }b>a\ast(1-\varepsilon)$$


Under certain conditions, the following result implies that a Cauchy subsequence exists. 

\begin{theorem}\label{th:main1}
Let $\{(X_n,M_n,\ast)\}_{n\in\mathbb{N}}$ be a sequence of nonempty compact fuzzy metric spaces satisfying the following conditions:
\begin{enumerate}
\item The $t$-norm $\ast$ satisfies property (TN1).
\item There exists a nondecreasing left-continuous function $C:(0,+\infty)\rightarrow (0,1]$, with 
$$0<C(s)\leq diam_s(X_n), \quad\text{for all $s>0$ and for all $n\in \mathbb{N}$}.$$
\item For every  $t>0$ and $0<\varepsilon<1$, there exists a natural number $N(\varepsilon,t)$ such that $Cov(X_n,\varepsilon,t)\leq N(\varepsilon,t)$, for all $n\in \mathbb{N}.$  

\item For every $t>0$ and  $0<\varepsilon<1$, consider $N=N(\varepsilon,t)$ as in Condition (3).  Assume that there exists a sequence  $\{\{x_i^{n}\}_{i=1}^{N}\}_{n\in\N}$ with  $\{x_i^{n}\}_{i=1}^N$ a $(t,\varepsilon)$-net in $X_{n}$, such that for every $n,m\in\N$ the following property is satisfied:  for all $s>t$ and each $i,j\in\{1,2,\ldots,N\}$, if $M_{n}(x_i^{n},x_j^{n},s)<M_{m}(x_i^{m},x_j^{m},s)$ then 

\[
\dfrac{M_{n}(x_i^{n},x_j^{n},s)}{M_{m}(x_i^{m},x_j^{m},s)\ast(1-\varepsilon)} \geq \dfrac{M_{n}(x_i^{n},x_j^{n},t)}{M_{m}(x_i^{m},x_j^{m},t)\ast(1-\varepsilon)}\, .
\]
\end{enumerate}
Then, there exists a subsequence $\{(X_{n_k},M_{n_k},\ast)\}_{k\in\mathbb{N}}$, such that 
$$M_{GH}(X_{n_j}, X_{n_k},t)> (1-\varepsilon)\ast(1-\varepsilon),$$ for all $j,k\in\N$. 
\end{theorem}

\begin{proof}
Let  $t>0$ and $0<\varepsilon<1$ be given. Take, for every $n\in \N$, the $(t,\varepsilon)$-net in $X_n$, say $\{x_i^n\}_{i=1}^N$, given by (4).
Since, for every $i,j\in\{1,2,\ldots, N\}$, condition (2)  tells us that $M_n(x_i^n,x_j^n,t)\geq C(t)$, we have 
\[
1\leq \dfrac{M_n(x_i^n,x_j^n,t)}{C(t)\ast \varepsilon}\leq \dfrac{1}{C(t)\ast \varepsilon}\, .
\]
The previous inequalities do not  depend on of our choice of $n$,  so the matrix of the integer parts
\[
A_{ij}^n =\left[ \dfrac{M_n(x_i^n,x_j^n,t)}{C(t)\ast \varepsilon}\right]
\]
can take only a finite number of values. Thus, one of the rows $A_{ij}^n$ is repeated infinitely often. Therefore,  we can find a subsequence  $\{X_{n_k}\}_{k\in\mathbb{N}}$ such that 
$$A_{ij}^{n_k}=A_{ij}^{n_p}\text{ for all $k,p$}.$$
Renaming the subsequence as $\{X_n\}_{n\in\N}$, we just have shown that, for every $n,m$ and every $i,j$,
\[
\left[ \dfrac{M_n(x_i^n,x_j^n,t)}{C(t)\ast \varepsilon}\right]=\left[ \dfrac{M_m(x_i^m,x_j^m,t)}{C(t)\ast \varepsilon}\right], 
\]
so that 
\[
\Big| \dfrac{M_n(x_i^n,x_j^n,t)}{C(t)\ast \varepsilon} -  \dfrac{M_m(x_i^m,x_j^m,t)}{C(t)\ast \varepsilon}\Big| <1
\]
or equivalently,
\[
\Big| M_n(x_i^n,x_j^n,t) -  M_m(x_i^m,x_j^m,t)\Big| <C(t)\ast \varepsilon.
\]
By applying Lemma~\ref{lem:ab}  we get, for every $i$ and $j$,
\begin{equation}\label{eq:cond2}
\begin{array}{lll}
(a) & M_n(x_i^n,x_j^n,t) & > \, \, M_m(x_i^m,x_j^m,t)\ast(1-\varepsilon), \\
 (b) &  M_m(x_i^m,x_j^m,t)&> \, \, M_n(x_i^n,x_j^n,t)\ast(1-\varepsilon). 
\end{array}
\end{equation}
Now we take $i,j\in\{1,2,\ldots,N\}$ and we fix $n,m\in\N$.
We will prove that condition (2) of Proposition~\ref{admisible} is satisfied as well.
Fix $s>t$. 

If  $M_n(x_i^n,x_j^n,s)\geq M_m(x_i^m,x_j^m,s)$, then $M_n(x_i^n,x_j^n,s)> M_m(x_i^m,x_j^m,s)\ast(1-\varepsilon)$.
Moreover, if  $M_n(x_i^n,x_j^n,s)<M_m(x_i^m,x_j^m,s)$, then by applying the condition (4) and Equation~\eqref{eq:cond2}(a), we get  

\begin{equation}\label{eq:quocient}
\dfrac{M_n(x_i^n,x_j^n,s)}{M_m(x_i^m,x_j^m,s)\ast(1-\varepsilon)} \geq \dfrac{M_n(x_i^n,x_j^n,t)}{M_m(x_i^m,x_j^m,t)\ast(1-\varepsilon)}>1 
\end{equation}

which means that 
$M_n(x_i^n,x_j^n,s)> M_m(x_i^m,x_j^m,s)\ast(1-\varepsilon)$.

Therefore, we have  
$$M_n(x_i^n,x_j^n,s)\geq M_m(x_i^m,x_j^m,s)\ast(1-\varepsilon),$$
for all $s\geq t$.  A similar argument to the previous one shows  
$$M_m(x_i^m,x_j^m,s)\geq M_n(x_i^n,x_j^n,s)\ast(1-\varepsilon),$$
for all $s\geq t$. 

Now, by Proposition~\ref{admisible}, we get an admissible fuzzy metric $M$ on $X_n\sqcup X_m$ such that
$$H_M(X_n,X_m,t)>(1-\varepsilon)\ast(1-\varepsilon),$$
which yields $M_{GH}(X_n,X_m,t)>(1-\varepsilon)\ast(1-\varepsilon)$. This completes the proof. 
\end{proof}

According to the previous result and Remark~\ref{rem:Cauchy}, we obtain 

\begin{corollary}
If $\{(X_n, M_n, \ast)\}_{n\in\N}$ is a sequence of nonempty compact fuzzy metric spaces which satisfies the conditions of the previous theorem, then  it has a Cauchy subsequence. 
\end{corollary}

We can give a slightly improved version by using the product $t$-norm. Notice that the product $t$-norm satisfies property (TN1). Moreover, we can reformulate condition (4), since the inequality~(\ref{eq:quocient}) in the proof above is equivalent to 
\[
\dfrac{M_n(x_i^n,x_j^n,s)}{M_m(x_i^m,x_j^m,s)} \geq \dfrac{M_n(x_i^n,x_j^n,t)}{M_m(x_i^m,x_j^m,t)}>(1-\varepsilon). 
\]
With these remarks in mind we can show
\begin{theorem}\label{th:main2}
Let $\{(X_n,M_n,\cdot)\}_{n\in\mathbb{N}}$ be a sequence of nonempty compact fuzzy metric spaces satisfying the following conditions:
\begin{enumerate}
\item There exists a nondecreasing left-continuous function $C:(0,+\infty)\rightarrow (0,1]$, with 
$$0<C(s)\leq diam_s(X_n), \quad\text{for all $s>0$ and for all $n\in \mathbb{N}$}.$$
\item For every  $t>0$ and $0<\varepsilon<1$, there exists a natural number $N(\varepsilon,t)$ such that, for all $n\in \mathbb{N}$, the cover-number 
$$Cov(X_n,t,\varepsilon)\leq N(\varepsilon,t).$$

 \item For every $t>0$ and $0<\varepsilon<1$, consider $N=N(\varepsilon,t)$. Assume that there exists a sequence  $\{\{x_i^{n}\}_{i=1}^{N}\}_{n\in\N}$ with  $\{x_i^{n}\}_{i=1}^N$ a $(t,\varepsilon)$-net in $X_{n}$, such that for every $n,m\in\N$ the following property is satisfied:  for all $s>t$ and each $i,j\in\{1,2,\ldots,N\}$, if $M_{n}(x_i^{n},x_j^{n},s)<M_{m}(x_i^{m},x_j^{m},s)$ then 
$$
\dfrac{M_n(x_i^n,x_j^n,s)}{M_m(x_i^m,x_j^m,s)} \geq \dfrac{M_n(x_i^n,x_j^n,t)}{M_m(x_i^m,x_j^m,t)}.
$$
\end{enumerate}
Then,
given $t>0$ and $0<\varepsilon<1$, there exists a subsequence $\{(X_{n_k},M_{n_k},\cdot)\}_{k\in\mathbb{N}}$, such that 
$$M_{GH}(X_{n_j}, X_{n_k},t)> (1-\varepsilon)\cdot (1-\varepsilon),$$ for all $j,k\in\N$. 
Therefore, there exists a Cauchy subsequence of $\{(X_n,M_n,\cdot)\}_{n\in\mathbb{N}}$.
\end{theorem}

We will see in the next example that condition (3) in the above theorem cannot be dropped out.

\begin{example}\label{ex:NoCauchyseq}
For every $n\in\N$, let $X_n=\{x_1^n,x_2^n\}$ be a metric space endowed with the metric $d_n$ defined by:
$$d_n(x_1^n,x_2^n)=d_n(x_2^n,x_1^n)=n \quad\text{and}\quad  d_n(x_i^n,x_i^n)=0,\ i=1,2.$$

Consider now the sequence of compact fuzzy metric spaces $\{(X_n,M_n,\cdot)\}_{n\in \N}$
where 
$M_n(x_i^n,x_j^n,0)=0$ and, for all $n=1,2,\ldots$, $i,j=1,2$ and $t>0$
\begin{align*}
M_{n}(x_i^{n},x_j^{n},t)&=\begin{cases}
    \frac{1}{2} &\text{ if $t\leq d_n(x_i^{n},x_j^n)$}\\[1ex]
    1 &\text{ if $t> d_n(x_i^n,x_j^n)$}
    \end{cases}
    \qquad\text{if $n$ is even,}\\
    M_{n}(x_i^n,x_j^n,t)&=\begin{cases}
    \frac{1}{3} &\text{ if $t\leq d_n(x_i^n,x_j^n)$}\\[1ex]
    1 &\text{ if $t> d_n(x_i^n,x_j^n)$}
    \end{cases}
    \qquad\text{if $n$ is odd}.
\end{align*}

These fuzzy metrics are non-Archimedean (see Example~\ref{ex:FMdiscrete}).

This sequence satisfies condition (1) in the above theorem. Indeed,  $diam_s(X_n)=\frac{1}{2}$ for $n\in \N$ even and $diam_s(X_n)=\frac{1}{3}$ for $n\in \N$ odd. So we can take $C(s)=\frac{1}{3}$ for every $s>0$.
Condition (2) is also trivially satisfied, taking $N(\vr,t)=2$ for every $0<\vr<1$ and $t>0$.
Let us see that condition (3) fails to be true. For, 
consider $t=\frac{1}{2}<n=d_n(x_1^n,x_2^n)$ for every $n\in\N$, take any $m\in\N$, $n=2m$, and fix $n>s>m$. 
Then 
$$M_n(x_1^n,x_2^n,s)=\frac{1}{2}<1=M_m(x_1^m,x_2^m,s)$$
and 
$$
\dfrac{M_n(x_1^n,x_2^n,s)}{M_m(x_1^m,x_2^m,s)}=\frac{\frac{1}{2}}{1}  <\dfrac{\frac{1}{2}}{\frac{1}{2}}= \dfrac{M_n(x_1^n,x_j^n,t)}{M_m(x_i^m,x_j^m,t)}\quad\text{if $m$ is even},
$$
or
$$
\dfrac{M_n(x_1^n,x_2^n,s)}{M_m(x_1^m,x_2^m,s)}=\frac{\frac{1}{2}}{1}  <\dfrac{\frac{1}{2}}{\frac{1}{3}}= \dfrac{M_n(x_1^n,x_j^n,t)}{M_m(x_i^m,x_j^m,t)}\quad\text{if $m$ is odd}.
$$
Finally, let us show that the sequence of compact fuzzy metric spaces 
$\{(X_n,M_n,\cdot)\}_{n\in \N}$
does not have any Cauchy subsequence. Assume, on the contrary, there is a Cauchy subsequence which we rename  as $\{X_n\}_{n\in\N}$. Then, for $t=\frac{1}{2}$ and and $\varepsilon=\frac{1}{10}$, there exists $n_0$ such that 
$M_{GH}(X_n,X_m,\frac{1}{2})>1-\frac{1}{10}$, for all $n,m\geq n_0$.
Take $2m+1>2m\geq n_0$. Since the only $(\frac{1}{2},\frac{1}{10})$-net in each $X_n$ is $\{x_1^n,x_2^n\}$, Proposition~\ref{prop:necesaria} applies to get an admissible metric $M$ on $X_{2m}\sqcup X_{2m+1}$ such that
\begin{align*}
M\left(x_1^{2m},x_2^{2m},\frac{1}{2}\right)&\geq M\left( x_1^{2m+1},x_2^{2m+1},\frac{1}{2}\right)\cdot \left(\frac{9}{10}\right)^2,\\
M\left(x_1^{2m+1},x_2^{2m+1},\frac{1}{2}\right)&\geq M\left(x_1^{2m},x_2^{2m},\frac{1}{2}\right)\cdot \left(\frac{9}{10}\right)^2.\\
\end{align*}
This gives   
$\frac{1}{3}\geq \frac{1}{2}\cdot\left(\frac{9}{10}\right)^2$,
which is a contradiction.
\end{example}

Next, we will see that the condition $C(s)>0$ for all $s>0$ in Theorem~\ref{th:main2} cannot be relaxed either.
\begin{example}
 Consider the same spaces $(X_n,d_n)$, $n\in\N$, as in the example above, where the fuzzy metrics are now given by the standard fuzzy metric 
 $$M_{d_n}(x_i^n,x_j^n,t)=\dfrac{t}{t+d_n(x_i^n,x_j^n)}$$
 for all $t\geq 0$ and $n\in\N$. We will see that the sequence of compact non-Archimedean fuzzy metric spaces $\{X_n, M_{d_{n}}, \cdot\}$ has no Cauchy subsequences. 

 It is apparent that 
 $C(s)\leq\frac{s}{s+n}$
 for all $n\in\N$, giving $C(s)=0$ for all $s>0$.

 Assume that there exists a $M_{GH}$-Cauchy subsequence, say $(X_{n_k},M_{d_{n_k}},\cdot)$, $k\in\N$. By 
 \cite[Lemma 4.16]{MS:2015}, this gives also a $d_{GH}$-Cauchy sequence $\{(X_{n_k},d_{n_k})\}_{k\in\N}$ with 
 ${\rm diam}(X_{n_k})=n_k$. But that is not possible. Indeed,  \cite[Exercise 7.3.14]{BuBu} shows that,  for every metric spaces $X$ and $Y$, we have 
 $$d_{GH}(X,Y)\geq \frac{1}{2}\Big| {\rm diam}(X)-{\rm diam}(Y)\Big|.$$
\end{example}

For stationary metrics we can state the following result.  

\begin{theorem}\label{th:mainStat}
If a sequence $\{(X_n,M_n,\ast)\}_{n\in\mathbb{N}}$ of  nonempty compact fuzzy metric spaces sa\-tis\-fies the following three conditions:
\begin{enumerate}
\item The t-norm $\ast$ satisfies property (TN1) and  $(M_n,\ast)$ are stationary fuzzy metrics for all $n\in\N$.
\item There exists a constant $C>0$ with 
$$0<C \leq diam_s(X_n)=diam(X_n),\text{  for all $n\in \mathbb{N}$}.$$
\item For every $0<\varepsilon<1$, there exists a natural number $N(\varepsilon)$ such that, for all $n\in \mathbb{N}$, the cover-number 
$Cov(X_n,\varepsilon)\leq N(\varepsilon)$,
\end{enumerate}
then,
given  $0<\varepsilon<1$, there exists a subsequence $\{(X_{n_k},M_{n_k},\ast)\}_{k\in\mathbb{N}}$, such that 
$$M_{GH}(X_{n_j}, X_{n_k})> (1-\varepsilon)\ast (1-\varepsilon),$$ for all $j,k\in\N$. 
Therefore, there exists a Cauchy subsequence of $\{(X_n,M_n,\ast)\}_{n\in\mathbb{N}}$.
\end{theorem}
\begin{proof}
   Take  $0<\vr< 1$ and fix an $\vr$-net in each of the spaces $X_n$. Consider $N=N(\vr)\geq Cov(X_n,\vr)$, for all $n\in\N$ and  assume that each of these  $\vr$-nets is given by $\{x_k^n\}_{k=1}^{N}$, $n\in\N$. Proceed as in the proof of Theorem~\ref{th:main1} to get a subsequence, renamed as $\{X_n\}_{n\in\N}$, such that for all $i,j\in\{1,\ldots,N\}$ and every $n,m\in\N$ we have:
   \[
\Big| M_n(x_i^n,x_j^n,t) -  M_m(x_i^m,x_j^m,t)\Big| <C\ast \varepsilon.
\]
Now, apply Lemma~\ref{lem:ab} to get, for every $i$ and $j$,
\begin{equation*}\label{eq:cond3}
\begin{aligned}
(a)& \quad M_n(x_i^n,x_j^n,t)&> &\quad M_m(x_i^m,x_j^m,t)\ast(1-\varepsilon), \\
(b) & \quad M_m(x_i^m,x_j^m,t)&> &\quad M_n(x_i^n,x_j^n,t)\ast(1-\varepsilon),
\end{aligned}
\end{equation*}
which is true for every $t>0$, since all $(M_n,\ast)$ is a stationary fuzzy metrics.

Proposition~\ref{admisible} allows us to obtain  an admissible fuzzy metric $M$ on $X_n\sqcup X_m$ such that
$$H_M(X_n,X_m,t)>(1-\varepsilon)\ast(1-\varepsilon),$$
which yields $M_{GH}(X_n,X_m,t)>(1-\varepsilon)\ast(1-\varepsilon)$.
So, given $0<\varepsilon<1$, we have found a subsequence, say $\{X_n\}_{n\in\N}$, satisfying that, for all $n,m\in\N$,  
$$M_{GH}(X_n,X_m,t)>(1-\varepsilon)\ast(1-\varepsilon),$$
for all $t>0$,
so that, by Remark~\ref{rem:Cauchy},  there exists a Cauchy subsequence.
\end{proof}
%

We can now state a fuzzy version of the Gromov-Hausdorff compactness theorem. Recall that if $\mathcal{M}=\{(X\sb{j},M\sb{X\sb{j}},\ast) : j\in J\}$ is a set of compact
non-Archimedean fuzzy metric spaces, then we denote by $\mathcal{I}(\mathcal{M})$ the set of all isometry classes of elements of $\mathcal{M}$. $\mathcal{I}(\mathcal{M})$ is a metric space with the non-Archimedean Gromov-Hausdorff fuzzy metric.
\begin{theorem}\label{th:GH}
Let $\mathcal{M}$ be a family of nonempty compact non-Archimedean fuzzy metric spaces. Assume that the fuzzy metrics carry on a t-norm with the property (TN1) and
\begin{enumerate}
\item There exists a nondecreasing left-continuous function $C:(0,+\infty)\rightarrow (0,1]$, with 
$$0<C(s)\leq diam_s(X), \quad\text{for all $s>0$ and for all $X\in \mathcal{M}$}.$$
\item For every  $t>0$ and $0<\varepsilon<1$, there exists a number $N(t,\varepsilon)$ such that, for all $X\in \mathcal{M}$, $Cov(X,t,\varepsilon)\leq N(t,\varepsilon)$.


\item Every sequence in $\mathcal{M}$ satisfies Condition~(4) in Theorem~\ref{th:main1}.
\end{enumerate}
Then, $(\mathcal{I}(\mathcal{M}),M_{GH},\ast)$  is a compact space.
\end{theorem}
\begin{proof}
We know from Theorem~\ref{th:main1} that every sequence in $\mathcal{M}$  has a Cauchy subsequence which means that $\mathcal{I}(\mathcal{M})$ is a precompact space. In \cite{MS:2015} the authors have proved that this space is complete as well. Thus, $(\mathcal{I}(\mathcal{M}), M_{GH},\ast)$  is a compact space.
\end{proof}

\begin{remark}
It is worth noting that Condition (3) in the above theorem cannot be removed, since you can take a family $\mathcal{M}$ of compact non-Archimedean fuzzy metric spaces containing  the sequence of Example~\ref{ex:NoCauchyseq} and so $\mathcal{I}(\mathcal{M})$ cannot be a compact space.   
\end{remark}

Now, from our previous results we can deduce the classical Gromov-Hausdorff compactness theorem for metric spaces. Firstly we will need the following result.

\begin{theorem}\label{th:main3}
Let $\{(X_n,d_n)\}_{n\in\mathbb{N}}$ be a sequence of nonempty compact  metric spaces satisfying the following conditions:
\begin{enumerate}
\item There exists a constant $K>0$ such that 
$$ diam(X_n)\leq K, \quad\text { for all $n\in \mathbb{N}$}.$$
\item For every   $\varepsilon>0$, there exists a number $N(\varepsilon)$ such that, for all $n\in \mathbb{N}$, the cover-number 
$$Cov(X_n,\varepsilon)\leq N(\varepsilon).$$
\end{enumerate}
Then, the sequence $\{X_n,M_{d_n},\cdot)\}_{n\in\mathbb{N}}$ satisfies the conditions in Theorem~\ref{th:main2} and therefore it has a Cauchy subsequence.
\end{theorem}
\begin{proof}
Since $ diam(X_n)\leq K$, then, for every $s>0$ and for all $n\in\mathbb{N}$, we have that 
$diam_s(X_n)\geq \dfrac{s}{s+K}=C(s)>0$ which turns out to be a continuous and nondecreasing function.

On the other hand, given $0<\varepsilon<1$ and $t>0$, we have that 
$$Cov(X_n,\frac{\varepsilon t}{1-\varepsilon})\leq N(\frac{\varepsilon t}{1-\varepsilon})=N(\varepsilon,t)$$
and, consequently,  for all $n$, we have 
$$Cov(X_n,\varepsilon,t)\leq N(\varepsilon,t).$$

Finally, given $t>0$ and $0<\varepsilon<1$, $n,m\in\mathbb{N}$ and $N=N(\varepsilon,t)$, 
take any two $(t,\varepsilon)$-nets  $\{x_i^n\}_{i=1}^{N}$  and $\{x_i^m\}_{i=1}^{N}$ in $X_n$ and $X_m$,  respectively.
Let $s>t$ be given and assume that $M_n(x_i^n,x_j^n,s)<M_m(x_i^m,x_j^m,s)$. Then 
$
d_m(x_i^m,x_j^m)<d_n(x_i^n,x_j^n) 
$
and 
\[
\begin{split}
\dfrac{M_n(x_i^n,x_j^n,s)}{M_m(x_i^m,x_j^m,s)}&=\dfrac{s+d_m(x_i^m,x_j^m)}{s+d_n(x_i^n,x_j^n)} =1-\dfrac{d_n(x_i^n,x_j^n)-d_m(x_i^m,x_j^m)}{s+d_n(x_i^n,x_j^n)}\\[2ex]
&\geq 1-\dfrac{d_n(x_i^n,x_j^n)-d_m(x_i^m,x_j^m)}{t+d_n(x_i^n,x_j^n)}=\dfrac{t+d_m(x_i^m,x_j^m)}{t+d_n(x_i^n,x_j^n)}\\[2ex]
&=\dfrac{M_n(x_i^n,x_j^n,t)}{M_m(x_i^m,x_j^m,t)}.
\end{split}
\]
Thus, all three conditions in Theorem~\ref{th:main2} are fulfilled and we get the conclusion.
\end{proof}

\begin{theorem}\label{th:GHclassic}
Let $\mathcal{M}$ be a family of nonempty compact metric spaces satisfying
\begin{enumerate}
\item There exists a constant $K>0$ such that 
$$ diam(X)\leq K, \quad\text { for all $X\in \mathcal{M}$}.$$
\item For every   $\varepsilon>0$, there exists a number $N(\varepsilon)$ such that, for all $X\in \mathcal{M}$, the cover-number 
$$Cov(X,\varepsilon)\leq N(\varepsilon).$$
\end{enumerate}
Then, $\mathcal{I}(\mathcal{M})$ endowed with the Gromov-Hausdorff  metric is a compact space.
\end{theorem}
\begin{proof}
 We know that $(\mathcal{I}(\mathcal{M}),d_{GH})$ is a complete metric  space  \cite[Proposition 6]{pe}. Take any sequence $\{(X_n,d_n)\}_{n\in\N}$ in $\mathcal{M}$. By Theorem~\ref{th:main3}, we know that $\{(X_n,M_{d_n},\cdot)\}_{n\in\N}$ has a $M_{GH}$-Cauchy subsequence.
 Applying \cite[Lemma~4.16]{MS:2015} we have that it is also a $d_{GH}$-Cauchy subsequence, so that $\mathcal{I}(\mathcal{M})$ is a precompact space.
  Therefore,  $(\mathcal{I}(\mathcal{M}),d_{GH})$ is a  compact metric space.

\end{proof}
\section{Conclusion}

In this paper we present a fuzzy version of the Gromov's precompactness theorem which allows us to deduce the classical theorem for metric spaces by means of the standard fuzzy metric associated to a metric space. As expected, some additional conditions must be imposed to obtain the result, and we also show that they cannot be dropped out. Considering the applications of Gromov's \emph{classical} compactness theorem, our result opens the possibility of future applications in various fields of fuzzy mathematics.

 \end{document}